\documentclass[12pt]{report}
\usepackage{amssymb,latexsym,epsfig,mathrsfs,graphpap,graphics,amsmath,amssymb}

\begin{document}

\begin{flushleft}
 {\bf\Large { Nonuniform Periodic Wavelet Frames \parindent=0mm \vspace{.1in} on  Non-Archimedean  Fields}}

\parindent=0mm \vspace{.3in}

  {\bf{O. Ahmad$^{1}$, N. A. Sheikh$^{2}$ }}
\end{flushleft}

\parindent=0mm \vspace{.1in}
{{\it\small$^{1}${Department of Mathematics, National Institute of Technology, Srinagar-190006, Jammu and Kashmir, India. E-mail: $\text{siawoahmad@gmail.com}$}}

{{\it\small$^{2}${Department of Mathematics, National Institute of Technology, Srinagar-190006, Jammu and Kashmir, India. E-mail: $\text{neyaznit@yahoo.co.in}$}}

\parindent=0mm \vspace{.2in}
{\small
{\bf{Abstract}}: In real life application all signals are not obtained from  uniform shifts; so there is a natural question  regarding  analysis and decompositions of  these types of signals by a stable mathematical tool. Gabardo and Nashed  and  Gabardo and  Yu  filled this gap by  the concept of  nonuniform multiresolution analysis and  nonuniform wavelets based on the theory of spectral pairs  for which the associated translation set $\Lambda =\left\{ 0,r/N\right\}+2\,\mathbb Z$ is no longer a discrete subgroup of $\mathbb R$ but a spectrum associated with a certain one-dimensional spectral pair and the associated dilation is an even positive integer related to the given spectral pair. In this paper, we introduce a notion of nonuniform periodic wavelet frame on non-Archimedean field. Using Fourier transform technique and the unitary extension principle, we propose an approach  for the construction of nonuniform periodic wavelet frames on non-Archimedean fields.

\parindent=0mm \vspace{.2in}


\parindent=0mm \vspace{.2in}
 }

\parindent=0mm \vspace{.2in}

{\bf{1. Introduction}}

The notion of frames was first introduced by Duffin and Shaeffer \cite{DS} in connection with some deep problems in non-harmonic Fourier series, and more particularly with the question of determining when a family of exponentials $\left\{e^{i\alpha_{n}t}:n\in\mathbb Z\right\}$ is complete for $L^2[a, b]$. Frames did not generate much interest outside non- harmonic Fourier series until the seminal work by Daubechies et al .\cite{DGM}. After their pioneer work, the theory of frames began to be studied widely and deeply, particularly in the more specialized context of wavelet frames and Gabor frames. Frames provide a useful model to obtain signal decompositions in cases where redundancy, robustness, over-sampling, and irregular sampling ploy a role. Today, the theory of frames has become an interesting and fruitful field of mathematics with abundant applications in signal processing, image processing, harmonic analysis, Banach space theory, sampling theory, wireless sensor networks, optics, filter banks, quantum computing, and medicine. An important example about frame is wavelet frame, which is obtained by translating and dilating a finite family of functions. One of the most useful methods to construct wavelet frames is through the concept of unitary extension principle (UEP) introduced by Ron and Shen \cite{RS}  and were subsequently extended by Daubechies et al.\cite{DHRS} in the form of the Oblique Extension Principle (OEP). They give sufficient conditions for constructing tight and dual wavelet frames for any given refinable function $\phi(x)$ which generates a multiresolution analysis. Gabardo and Nashed \cite{GN1} and  Gabardo and  Yu \cite{GN2} notion of  nonuniform multiresolution analysis and  nonuniform wavelets based on the theory of spectral pairs  for which the associated translation set $\Lambda =\left\{ 0,r/N\right\}+2\,\mathbb Z$ is no longer a discrete subgroup of $\mathbb R$ but a spectrum associated with a certain one-dimensional spectral pair and the associated dilation is an even positive integer related to the given spectral pair. 

\parindent=8mm\vspace{.2in}
 
In recent years there has been a considerable interest in the problem of constructing periodic wavelet bases and frames in Hilbert spaces as most of the signals of practical interest are periodic in nature. Apart from signals that are inherently periodic, all signals resulting from experiments with a finite duration can in principle be modeled as periodic signals \cite{LP}. Since the setup of tight wavelet frames provides great flexibility in approximating and representing periodic functions. Using periodization techniques, Zhang \cite{zhang} constructed a dual pair of periodic wavelet frames for $L^2[0,1]$ under the assumption that the support of the wavelet function $\psi$ in the frequency domain  is contained in $[-\pi, -\varepsilon]\cup [\varepsilon, \pi],\, \varepsilon>0$. Zhang and Saito \cite{ZS} have constructed general periodic wavelet frames for $L^2[0, 1]$ using extension principles. Later on, Lu and Li \cite{LL} constructed periodic wavelet frames with dilation matrix. 

\parindent=8mm \vspace{.2in}
On the other hand, the past decade has also witnessed a tremendous interest in the problem of constructing wavelet bases and frames on various spaces other than $\mathbb R$. For example, R. L. Benedetto and J. J. Benedetto \cite{BB} developed a wavelet theory for local fields and related groups. They did not develop the multiresolution analysis (MRA) approach, their method is based on the theory of wavelet sets and only allows the construction of wavelet functions whose Fourier transforms are characteristic functions of some sets. Jiang et al.\cite{jiang} pointed out a method for constructing orthogonal wavelets on  local field $K$ with a constant generating sequence and derived necessary and sufficient conditions for a solution of the refinement equation to generate a multiresolution analysis of $L^2(K)$. Subsequently, tight wavelet frames on local fields of positive characteristic were constructed by Shah and Debnath \cite{SDtight} using extension
principles. Ahmad and Sheikh \cite{onuwf} introduced the notion of nonuniform wavelet frame on local fields and subsquently established a complete characterization of tight wavelet frames on local fields by virtue of two basic equations in the frequency domain and show how to construct an orthonormal wavelet basis for $L^2(\mathbb K)$. For more about frames, we refer to \cite{{chrisBOOK},{oJMP},{oJGP}}

\pagestyle{myheadings}

\parindent=8mm \vspace{.2in}
Drawing inspiration from the above work, our aim is to extend the notion of wavelet frames to  nonuniform periodic wavelet frames on non-Archimedean fields via extension principles. More precisely, we prove that under some mild conditions, the periodization of any nonuniform wavelet frame constructed by the unitary extension principle is a nonuniform periodic wavelet frame on non-Archimedean fields. 

\parindent=8mm \vspace{.2in}
The layout of this paper is  as follows. In Section 2, we discuss some preliminary facts about non-Archimedean fields  and also some results which are required in the subsequent sections. Sections 3 is devoted to discuss our  main results about nonuniform periodic wavelet frames.

\newpage
\parindent=0mm \vspace{.2in}
{\bf{2. Preliminaries and Nonuniform Periodic Wavelet System on  Non-Archimedean  Fields }}

\parindent=0mm \vspace{.1in}
A  non-Archimedean field $K$ is a  locally compact, non-discrete and totally disconnected field. If it is of characteristic zero, then  it is a field of $p$-adic numbers $\mathbb Q_p$ or its finite extension. If $K$ is of positive characteristic, then $K$ is a field of formal Laurent series over a finite field $GF(p^c)$. If $c =1$, it is a $p$-series field, while for $c\ne 1$, it is an algebraic extension of degree  $c$ of a $p$-series field. Let $K$ be a fixed non-Archimedean  field  with the ring of integers ${\mathfrak D}= \left\{x \in K: |x| \le 1\right\}$. Since $K^{+}$ is a locally compact Abelian group, we choose a Haar measure $dx$ for $K^{+}$. The  field $K$ is locally compact, non-trivial, totally disconnected and complete topological field endowed with non--Archimedean norm  $|\cdot|:K\to \mathbb R^+$ satisfying

\parindent=0mm \vspace{.1in}
(a) $|x|=0$ if and only if $x = 0;$

\parindent=0mm \vspace{.1in}
(b) $|x\,y|=|x||y|$ for all $x, y\in K$;

\parindent=0mm \vspace{.1in}
(c) $|x+y|\le \max \left\{ |x|, |y|\right\}$ for all $x, y\in K$.

\parindent=0mm \vspace{.1in}
Property (c) is called the ultrametric inequality. Let ${\mathfrak B}= \left\{x \in K: |x| < 1\right\}$ be the prime ideal of the ring of integers ${\mathfrak D}$ in $K$. Then, the residue space ${\mathfrak D}/{\mathfrak B}$ is isomorphic to a finite field $GF(q)$, where $q = p^{c}$ for some prime $p$ and $c\in\mathbb N$. Since  $K$ is totally disconnected and $\mathfrak B$ is both prime and principal ideal, so there exist a prime element $\mathfrak p$ of $K$ such that ${\mathfrak B}= \langle \mathfrak p \rangle=\mathfrak p {\mathfrak D}$. Let ${\mathfrak D}^*= {\mathfrak D}\setminus {\mathfrak B }=\left\{x\in K: |x|=1   \right\}$. Clearly,  ${\mathfrak D}^*$ is a group of units in $K^*$ and if $x\not=0$, then can write $x=\mathfrak p^n y, y\in {\mathfrak D}^*.$ Moreover, if ${\cal U}= \left\{a_m:m=0,1,\dots,q-1 \right\}$ denotes the fixed full set of coset representatives of ${\mathfrak B}$ in ${\mathfrak D}$, then every element $x\in K$ can be expressed uniquely  as $x=\sum_{\ell=k}^{\infty} c_\ell \,\mathfrak p^\ell $ with $c_\ell \in {\cal U}.$ Recall that ${\mathfrak B}$ is compact and open, so each  fractional ideal ${\mathfrak B}^k= \mathfrak p^k {\mathfrak D}=\left\{x \in K: |x| < q^{-k}\right\}$  is also compact and open and is a subgroup of $K^+$. We use the notation in Taibleson's book [13]. In the rest of this paper, we use the symbols $\mathbb N, \mathbb N_0$ and $\mathbb Z$ to denote the sets of natural, non-negative integers and integers, respectively.

\parindent=8mm \vspace{.1in}
 Let $\chi$ be a fixed character on $K^+$ that is trivial on ${\mathfrak D}$ but  non-trivial on  ${\mathfrak B}^{-1}$. Therefore, $\chi$ is constant on cosets of ${\mathfrak D}$ so if $y \in {\mathfrak B}^k$, then $\chi_y(x)=\chi(y,x), x\in K.$ Suppose that $\chi_u$ is any character on $K^+$, then the restriction $\chi_u|{\mathfrak D}$ is a character on ${\mathfrak D}$. Moreover, as characters on ${\mathfrak D}, \chi_u=\chi_v$ if and only if $u-v\in {\mathfrak D}$. Hence, if  $\left\{u(n): n\in\mathbb N_0\right\}$ is a complete list of distinct coset representative of ${\mathfrak D}$ in $K^+$, then, as it was proved in [13], the set  $\left\{\chi_{u(n)}: n\in\mathbb N_0\right\}$   of distinct characters on ${\mathfrak D}$ is a complete orthonormal system on ${\mathfrak D}$.

\parindent=8mm \vspace{.1in}
We now impose a natural order on the sequence $\{u(n)\}_{n=0}^\infty$. We have ${\mathfrak D}/ \mathfrak B \cong GF(q) $ where $GF(q)$ is a $c$-dimensional vector space over the field $GF(p)$. We choose a set $\left\{1=\zeta_0,\zeta_1,\zeta_2,\dots,\zeta_{c-1}\right\}\subset {\mathfrak D^*}$ such that span$\left\{\zeta_j\right\}_{j=0}^{c-1}\cong GF(q)$. For $n \in \mathbb N_0$ satisfying
$$0\leq n<q,~~n=a_0+a_1p+\dots+a_{c-1}p^{c-1},~~0\leq a_k<p,~~\text{and}~k=0,1,\dots,c-1,$$

\parindent=0mm \vspace{.1in}
we define
$$u(n)=\left(a_0+a_1\zeta_1+\dots+a_{c-1}\zeta_{c-1}\right){\mathfrak p}^{-1}.$$

\parindent=0mm \vspace{.1in}
Also, for $n=b_0+b_1q+b_2q^2+\dots+b_sq^s, ~n\in \mathbb N_{0},~0\leq b_k<q,k=0,1,2,\dots,s$, we set

$$u(n)=u(b_0)+u(b_1){\mathfrak p}^{-1}+\dots+u(b_s){\mathfrak p}^{-s}.\eqno$$

\parindent=0mm \vspace{.1in}
This defines $u(n)$ for all $n\in \mathbb N_{0}$. In general, it is not true that $u(m + n)=u(m)+u(n)$. But, if $r,k\in\mathbb N_{0}\; \text{and}\;0\le s<q^k$, then $u(rq^k+s)=u(r){\mathfrak p}^{-k}+u(s).$ Further, it is also easy to verify that $u(n)=0$ if and only if $n=0$ and $\{u(\ell)+u(k):k \in \mathbb N_0\}=\{u(k):k \in \mathbb N_0\}$ for a fixed $\ell \in \mathbb N_0.$ Hereafter we use the notation $\chi_n=\chi_{u(n)}, \, n\ge 0$.

\parindent=8mm \vspace{.1in}
Let the local field $K$ be of characteristic $p>0$ and $\zeta_0,\zeta_1,\zeta_2,\dots,\zeta_{c-1}$ be as above. We define a character $\chi$ on $K$ as follows:
$$\chi(\zeta_\mu {\mathfrak p}^{-j})= \left\{
\begin{array}{lcl}
\exp(2\pi i/p),&&\mu=0\;\text{and}\;j=1,\\
1,&&\mu=1,\dots,c-1\;\text{or}\;j \neq 1.
\end{array}
\right. \eqno(2.1)$$
The Fourier transform of $f \in L^1(K)$ is denoted by $\hat f(\xi)$ and defined  by
\begin{align*}
{\mathscr F}\big\{f(x)\big\}=\hat f(\xi)=\int_K f(x)\overline{ \chi_\xi(x)}\,dx.\tag{2.2}
\end{align*}
It is noted that
$$\hat f(\xi)= \displaystyle \int_K f(x)\,\overline{ \chi_\xi(x)}dx= \displaystyle \int_K f(x)\chi(-\xi x)\,dx.\eqno(2.3)$$

\parindent=0mm \vspace{.1in}
The properties of Fourier transforms on non-Archimedean field $K$ are much similar to those of on the classical field $\mathbb R$. In fact, the Fourier transform on non-Archimedean fields of positive characteristic have the following properties:
\begin{itemize}
  \item The map $f\to \hat f$ is a bounded linear transformation of $L^1(\mathbb K)$ into $L^\infty(\mathbb K)$, and $\big\|\hat f\big\|_{\infty}\le \big\|f\big\|_{1}$.
  \item If $f\in L^1(\mathbb K)$, then $\hat f$ is uniformly continuous.
  \item If $f\in L^1(\mathbb K)\cap L^2(\mathbb K)$, then $\big\|\hat f\big\|_{2}=\big\|f\big\|_{2}$.
\end{itemize}

\parindent=0mm \vspace{.1in}
The Fourier transform of a function $f\in L^2(\mathbb K)$ is defined by
\begin{align*}
\hat f(\xi)= \lim_{k\to \infty} \hat f_{k}(\xi)=\lim_{k\to \infty}\int_{|x|\le q^{k}} f(x)\overline{ \chi_\xi(x)}\,dx,\tag{2.4}
\end{align*}

\parindent=0mm \vspace{.0in}
where  $f_{k}=f\,\Phi_{-k}$ and $\Phi_{k}$ is  the characteristic function of ${\mathfrak B}^{k}$.  Furthermore, if $f\in L^2(\mathfrak D)$, then we define the Fourier coefficients of $f$ as
\begin{align*}
\hat f\big(u(n)\big)=\int_{\mathfrak D} f(x) \overline{ \chi_{u(n)}(x)}\,dx.\tag{2.5}
\end{align*}

\parindent=0mm \vspace{.0in}
The series $\displaystyle\sum_{n\in \mathbb N_{0}} \hat f\big( u(n)\big) \chi_{u(n)}(x)$ is called the Fourier series of $f$. From the standard $L^2$-theory for compact Abelian groups, we conclude that the Fourier series of $f$ converges to $f$ in $L^2(\mathfrak D)$ and Parseval's identity holds:
\begin{align*}
\big\|f\big\|^2_{2}=\int_{\mathfrak D}\big|f(x)\big|^2 dx= \sum_{n\in \mathbb N_{0}} \left| \hat f\big(u(n)\big)\right|^2.\tag{2.6}
\end{align*}

\parindent=0mm \vspace{.1in}
 We also denote the test function space on $K$ by $\Omega(K)$, that is, each function $f$ in $\Omega(K)$ is a finite linear combination of functions of the form ${\bf 1}_k(x-h), h\in K, k\in\mathbb Z$, where ${\bf 1}_k$ is the characteristic function of ${\mathfrak B}^k$. This class of functions can also be described in the following way. A function $g\in \Omega(K)$ if and only if there exist integers $k,\ell$ such that $g$ is constant on cosets of ${\mathfrak B}^{k}$ and is supported on ${\mathfrak B}^{\ell}$. It follows that $\Omega$ is closed under Fourier transform and is an algebra of continuous functions with compact support, which is dense in ${\cal C}_{0}(K)$ as well as in $L^p(K),  1\le p <\infty$. For more details we refer to \cite{{RV},{table}}.
 
\parindent=8mm \vspace{.1in}
 For an integer $N \ge 1$ and an odd integer $r$ with $1\leq r \leq qN-1$ such that $r$ and $N$ are relatively prime, we define $$\Lambda = \left\{ 0, \dfrac{u(r)}{N}\right\}+{\mathcal Z}.$$

 \parindent=0mm \vspace{.1in}
where ${\mathcal Z}=\left\{ u(n): n\in \mathbb N_{0}\right\}$. It is easy to verify that $\Lambda$ is not a group on non-Archimedean field $K$, but is the union of ${\mathcal Z}$ and a translate of ${\mathcal Z}.$

%
%
%
%

\parindent=8mm \vspace{.1in}
As in the standard scheme, one expects the existence of $qN -1$ number of functions so that their translation by elements of $\Lambda$ and dilations by the integral powers of ${\mathfrak p^{-1}}N$ form an orthonormal basis for $L^2(\mathbb K)$.

\parindent=8mm \vspace{.1in}

For $j\in\mathbb N_{0}$, let ${\cal N}_{j}$ denote a full collection of coset representatives of $\Lambda/(qN)^{j}\Lambda$, i.e.,

$${\cal N}_{j}=\left\{0,1,2,\dots, (qN)^{j}-1\right\}, \quad j\ge 0.$$

\parindent=0mm \vspace{.1in}
Then, $\Lambda=\bigcup _{n\in {\cal N}_{j}} \left(n+(qN)^{j}\Lambda\right),$ and for any distinct $n_{1}, n_{2}\in {\cal N}_{j}$, we have $\left(n_{1}+(qN)^{j}\Lambda\right)\cap \left(n_{2}+(qN)^{j}\Lambda\right)=\emptyset.$ Thus, every non-negative integer $k$ can uniquely be written as
$k=r(qN)^{j}+s$, where $r\in \Lambda, s\in {\cal N}_{j}$. Further, a bounded function $W : K\to K$  is said to be a radial decreasing $L^1$-majorant of $f(x)\in L^2(\mathbb K)$ if $|f(x)|\le W(x),\,W\in L^{1}(\mathbb K),$ and  $W(0)<\infty.$

\parindent=8mm \vspace{.1in}
Let $a$ and $b$ be any two fixed elements in $K$. Then, for any  prime  $\mathfrak p$ and $m,n\in\mathbb N_{0}$,let $D_{\mathfrak p},T_{u(n) a}$ and $E_{u(m) b}$   be the unitary operators acting on $f\in L^2(\mathbb K)$ defined by :
\begin{align*}
& T_{u(n) a}f(x)=f\big(x-u(n) a\big), \qquad ~\text{(Translation)}\\
&E_{u(m) b}f(x)=\chi\big(u(m) b x\big) f(x),~\quad \text{(Modulation)}\\
&D_{\mathfrak p}f(x)=\sqrt {qN} f\left(\mathfrak p^{-1}N x\right), ~\quad\quad \text{(Dilation)}.
\end{align*}

\parindent=0mm \vspace{.1in}
For given $\Psi := \left\{\psi_1,\dots, \psi_{qN-1}\right\}\subset  L^2(K)$, define the  nonuniform wavelet system
$${\cal W}({\Psi},\lambda)=\Big\{ \psi_{\ell,j,\lambda}:=(qN)^{j/2}\psi_\ell\big(({\mathfrak p}^{-1}N)^j x-\lambda\big),\, j\in \mathbb Z,\lambda\in \Lambda, 1 \le \ell \le qN-1 \Big\}.\eqno(2.7)$$
The nonuniform wavelet system ${\cal W}({\Psi},\lambda)$ is called a  {\it non uniform wavelet frame}, if there exist positive numbers $0 < A \le B < \infty$ such that for all $f\in  L^2(\mathbb K)$

$$A\big\|f \big\|^2_{2} \le \sum_{\ell=1}^{qN-1}\sum_{j\in\mathbb Z}\sum_{\lambda \in \Lambda} \left|\big\langle f, \psi_{\ell,j, \lambda}\big\rangle\right|^2 \le B \big\|f\big\|^2_{2}.\eqno(2.8)$$

\parindent=0mm \vspace{.1in}
The largest $A$ and the smallest $B$ for which $(2.8)$ holds are called {\it non uniform  wavelet frame bounds}. A wavelet frame is a {\it tight non uniform  wavelet frame} if $A$ and $B$ are chosen such that $A = B$ and then the generators $\left\{\psi_1, \psi_2,\dots,\psi_{qN-1}\right\}$ are often referred as {\it tight non uniform  framelets}. If only the right-hand inequality in (2.8) holds, then ${\cal W}({\Psi},\lambda)$ is called a {\it Bessel sequence}.





\parindent=8mm \vspace{.1in}
Next, we give a brief account of the MRA-based wavelet frames generated by the wavelet masks on non-Archimedean local fields. Following the unitary extension principle, one often starts with a refinable function or even with a refinement mask to construct desired wavelet frames. A  function $\varphi\in L^2(\mathbb K)$ is called a {\it nonuniform refinable function}, if it satisfies an equation of the type
$$\varphi(x)= \sqrt{qN}\sum_{\lambda \in \Lambda} a_\lambda\,\varphi\big(({\mathfrak p}^{-1}N)x- \lambda\big),\eqno(2.9)$$

\parindent=0mm \vspace{.0in}
where $ \big\{a_\lambda: \lambda \in \Lambda \big\} \in l^2(\mathbb N_0).$  In the frequency domain,  equation (2.9) can be written as
$$ \hat\varphi\big( \xi \big)=  m_0\left(\dfrac{{\mathfrak  p } \xi}{N}\right)\hat \varphi\left(\dfrac{{\mathfrak  p } \xi}{N}\right),\eqno(2.10) $$

\parindent=0mm \vspace{.1in}
where $$m_0(\xi)=\dfrac{1}{\sqrt{qN}}\displaystyle\sum_{\lambda \in \Lambda }a_\lambda\, \overline{\chi_\lambda(\xi)},\eqno(2.11)$$

\parindent=0mm \vspace{.1in}
is an integral periodic  function in $L^2({\mathfrak  D })$ and is often called the {\it refinement mask} of $\varphi$.  Observe that $\chi_k(0)=\hat\phi(0)=1$. By letting $\xi=0$ in equations (2.10) and (2.11), we obtain  $\displaystyle\sum_{\lambda \in \Lambda} a_{\lambda}=1$. Further, it is proved in \cite{SDtight} that a function $\varphi\in L^2(\mathbb K)$ generates an MRA in $L^2(\mathbb K)$ if and only if
$$ \sum_{\lambda \in \Lambda}\left|\hat \varphi\big(\xi+\lambda \big)\right|^2=1,~\text{for}~a.e.~\xi \in \mathfrak D,~~ \text{and}~~ \hat \varphi(0)= \lim_{\xi\to 0} \hat \varphi(\xi)=1,\quad \xi\in \mathbb K.\eqno(2.12)$$

\parindent=0mm \vspace{.0in}
Suppose $\Psi=\left\{\psi_1, \dots , \psi_{qN-1}\right\}$ is a set of MRA functions derived from

$$\hat \psi_\ell(\xi)=m_\ell\left(\dfrac{{\mathfrak  p } \xi}{N}\right){\hat \varphi}\left(\dfrac{{\mathfrak  p } \xi}{N}\right),\eqno(2.13) $$

where
$$m_\ell(\xi)=\dfrac{1}{\sqrt {qN}}\displaystyle\sum_{\lambda \in \Lambda}a_\lambda^\ell\; \overline{\chi_\lambda(\xi)},\quad  1 \le \ell \le qN-1 \eqno(2.14) $$

\parindent=0mm \vspace{.1in}
are the integral periodic  functions in $L^2({\mathfrak  D })$ and are called the {\it nonuniform framelet symbols} or {\it nonuniform wavelet masks}. With  $m_\ell(\xi), 0 \le \ell \le qN-1$ as the wavelet masks, we formulate the matrix ${\cal M}(\xi)$ as

$${\mathcal M}(\xi)=\left(\begin{array}{cccc}
m_0(\xi)&m_0\big(\xi+ {\mathfrak p}u(1)\big)&\dots&m_0\big(\xi+ {\mathfrak p}u(s-1)\big)\\
m_1(\xi)&m_1\big(\xi+ {\mathfrak p}u(1)\big)&\dots&m_1\big(\xi+ {\mathfrak p}u(s-1)\big)\\
\vdots&\vdots&\ddots&\vdots\\
m_{qN-1}(\xi)&m_{qN-1}\big(\xi+ {\mathfrak p}u(1)\big)&\dots&m_{qN-1}\big(\xi+ {\mathfrak p}u(s-1)\big)\\
\end{array}\right).\eqno(2.15)$$

\parindent=0mm \vspace{.1in}
The so-called unitary extension principle (UEP) provides a sufficient condition on $\Psi=\left\{\psi_1,\psi_2. \dots , \psi_{qN-1}\right\}$ such that the  nonuniform wavelet system ${\mathcal W}({\Psi},\lambda)$ given by (2.7) constitutes a tight frame for $L^2(\mathbb K)$. It is well known that in order to apply the UEP to derive wavelet tight frame from a given refinable function, the corresponding refinement mask must satisfy

$$\sum_{\lambda=0}^{qN-1}\left|m_{0}\big(\xi+qN\lambda\big)\right|^{2}\le 1,\quad \xi\in \mathbb K.\eqno(2.16)$$

\parindent=8mm \vspace{.1in}
In this connection, Shah and Debnath \cite{SDtight} gave an explicit construction scheme for the construction of tight wavelet frames on local fields of positive characteristic using  unitary extension principles. The following is the fundamental tool they gave to construct tight wavelet frames on local fields.

\parindent=0mm \vspace{.2in}
{\bf{Theorem 2.1.}}  Let  $\varphi(x)$ be a compactly supported refinable function and $\hat\phi(0)=1$. Then, the nonuniform wavelet system ${\mathcal W}({\Psi},\lambda)$ given by (2.7) constitutes a normalized tight wavelet frame in $L^2(K)$ provided the matrix ${\mathcal M}(\xi)$ as defined in (2.15) satisfies 

$${\mathcal M}(\xi){\mathcal M^*}(\xi)=I_q,\quad \text{for}~a.e.~\xi \in \sigma (V_0)\eqno(2.17)$$

\parindent=0mm \vspace{.1in}
 where $\sigma (V_0):=\Big\{ \xi\in {\mathfrak D}:  \displaystyle\sum_{\lambda\in \Lambda}|\hat\varphi\big(\xi+\lambda\big)|^2\not=0  \Big\}.$  

\parindent=8mm \vspace{.2in}






\parindent=0mm \vspace{.1in}

{\bf{3. Nonuniform Periodic Wavelet Frames on Non-Archimedean Fields }}

\parindent=0mm \vspace{.2in}
In this section, we present an approach for constructing nonuniform periodic wavelet frames on non-Archimedean fields by virtue of unitary extension principle (UEP).

\parindent=0mm \vspace{.1in}
For any $f \in L^1(\mathbb K)$, we define  the periodic version of $f$ as

$$f^{\text{per}}=\sum_{\lambda \in \Lambda}f\big(x+\lambda u(N)\big).$$

\parindent=0mm \vspace{.1in}
It is easy to see that $f^{\text{per}}$ is well defined and it is $\mathbb N$-periodic local integrable function. With same dilation  and translation operators as defined in Section 2, we define the {\it nonuniform periodic wavelet system} as
$${\cal W}^{\text{per}}({\Psi},\lambda):=\Big\{\varphi^{\text{per}},\psi_{\ell, j, \lambda}^{\text{per}}:1\le \ell \le qN-1,j\in\mathbb N_{0},\lambda \in {\cal N}_{j}\Big\}.\eqno(3.1) $$

\parindent=8mm \vspace{.1in}

\parindent=8mm \vspace{.2in}
In order to establish the main result of this section, we first state and prove the following lemmas.

\parindent=0mm \vspace{.1in}

{\bf{Lemma 3.1.}}  Suppose that the nonuniform  periodic wavelet system ${\cal W}^{\text{per}}({\Psi},\lambda)$ is defined by (3.1). Then, for any periodic function $f $ and given $\varepsilon>0$, there exists a positive integer $J\in \mathbb N$, such that

$$(1-\varepsilon)\big\|f\big\|_2^2 \le \sum_{\lambda =0}^{(qN)^j -1}\left|\left\langle f, \varphi_{j,\lambda}^{\text{per}}\right\rangle\right|^2 \le  (1+\varepsilon)\big\|f\big\|_2^2, \quad \text{for all}~ j \ge J. \eqno(3.2)$$ 

\parindent=0mm \vspace{.1in}

{\bf{Proof.}} Let $\Gamma$ denotes the support of the Fourier coefficients $\left\{\hat f\left( u\left(\dfrac{r}{N}\right)\right)\right\}_{r \in \mathbb N_{0}}$. Then, we have

$$f(x)=\sum_{r\in \Gamma}\hat f\left( u\left(\dfrac{r}{N}\right)\right)\overline{ \chi_{r/N}(x)}.$$

\parindent=0mm \vspace{.1in}
Assume that
$$\varphi_{j,\lambda}^{\text{per}}(x)=\sum_{r\in \mathbb N_0}\hat\varphi_{j,\lambda}^{\text{per}}\left(u\left(\dfrac{r}{N}\right)\right)\overline{ \chi_{r/N}(x)}.\eqno(3.3)$$

\parindent=0mm \vspace{.1in}
The Fourier coefficients in the above  expression can be written as
$$\hat\varphi_{j,\lambda}^{\text{per}}\left(u\left(\dfrac{r}{N}\right)\right)=(qN)^{j/2}\hat\varphi\left(({\mathfrak p}^{-1}N)^{-j} u\left(\dfrac{r}{N}\right)\right)\,\overline{\chi_\lambda \left(({\mathfrak p}^{-1}N)^{-j} u\left(\dfrac{r}{N}\right)\right)}.\eqno(3.4)$$

\parindent=0mm \vspace{.1in}

Parseval's formula of the above Fourier series gives

$$\begin{array}{rcl}
\displaystyle\sum_{\lambda =0}^{(qN)^j -1}\left|\left\langle f, \varphi_{j,\lambda}^{\text{per}}\right\rangle\right|^2&=&\displaystyle\sum_{\lambda =0}^{(qN)^j -1}\Big|\sum_{r \in \Gamma}\hat f\big( u(r)\big)\overline{\hat\varphi_{j,\lambda}^{\text{per}}\big(u(r)\big)}\Big|^2\\\\
&=&\displaystyle\sum_{\lambda =0}^{(qN)^j -1}\Big|\sum_{r\in \Gamma}\hat f\big( u(r)\big)(qN)^{j/2}\overline{\hat\varphi\big(({{\mathfrak p}^{-1}N})^{-j} u(r)\big)}\chi_\lambda\big(({\mathfrak p}^{-1}N)^{-j}u(r)\big)\Big|^2\\\\
&=&\displaystyle\sum_{\lambda =0}^{(qN)^j -1} \Big|\sum_{r\in \Gamma} d_r\left(\hat f, \hat\varphi\right)\chi_\lambda\big(({\mathfrak p}^{-1}N)^{-j}u(r)\big)\Big|^2,
\end{array}$$

\parindent=0mm \vspace{.1in}
where $d_r\left(\hat f, \hat\varphi\right)=(qN)^{j/2} \hat f\big( u(r)\big)\overline{\hat\varphi\big(({\mathfrak p}^{-1}N)^{-j} u(r)\big)}.$ Since $|Gamma$ is a finite set, there exists a positive number $N$ such that $E \subseteq D(N):=\big\{\lambda \in \Lambda:\, |\lambda| \le N\big\}.$ Hence, there exists $J_1\ge 0$, such that for all $j \ge J_1$, the elements of $D(N)$ lie in different cosets of $\Lambda/(qN)^{j}\Lambda$. Thus, the cardinality of $\Gamma\cap (\lambda+(qN)^{j}\Lambda)$ is at most one for each $j\ge J_{1}, \lambda\in {\cal N}_{j}$. Consequently, we have
$$\begin{array}{rcl}
\displaystyle\sum_{\lambda =0}^{(qN)^j -1}\left|\left\langle f, \varphi_{j,\lambda}^{\text{per}}\right\rangle\right|^2&=&\displaystyle\sum_{\lambda =0}^{(qN)^j -1}\sum_{r\in \Gamma} d_r\left(\hat f, \hat\varphi\right)\chi_\lambda\big(({\mathfrak p}^{-1} N)^{-j}u(r)\big)\sum_{s\in \Gamma}\overline{d_s\left(\hat f, \hat\varphi\right)\chi_\lambda\big(({\mathfrak p}^{-1}N)^{-j}u(s)\big)}
\\\\
&=&\displaystyle\sum_{r\in \Gamma}\sum_{s\in \Gamma}d_r\left(\hat f, \hat\varphi\right)\overline{d_r\left(\hat f, \hat\varphi\right)}\sum_{\lambda =0}^{(qN)^j -1}\,\chi_\lambda\big(({\mathfrak p}^{-1}N)^{-j}u(r-s)\big)\\\\
&=&(qN)^{j}\displaystyle\sum_{s\in \Gamma}\Big|d_s\left(\hat f, \hat\varphi\right)\Big|^2\\\\
&=&\displaystyle\sum_{s\in \Gamma}\Big|(qN)^{j/2} \hat f\big( u(s)\big)\overline{\hat\varphi\big(({\mathfrak p}^{-1}N)^{-j} u(s)\big)}\Big|^2.
\end{array}$$

\parindent=0mm \vspace{.0in}
Since $\hat \varphi(0)=\lim_{\xi \to 0}\hat \varphi(\xi)=1$, therefore there exists a non-negative integer $J_2$, such that

$$(1-\varepsilon)\le \left|\hat\varphi\big(({\mathfrak p}^{-1}N)^{-j}u(s)\big)\right|^2 \le (1+\varepsilon), \quad \text{for all}~j\ge J_{2}. $$

\parindent=0mm \vspace{.1in}
Let $J=\max\left\{J_1, J_2\right\}$, then with this choice of $j \ge J$, we obtain

$$(1-\varepsilon)\sum_{s\in \Gamma}\left|\hat f\big( u(s)\big)\right|^2\le \sum_{\lambda =0}^{(qN)^j -1}\left|\left\langle f, \varphi_{j,\lambda}^{\text{per}}\right\rangle\right|^2\le(1+\varepsilon)\sum_{s\in \Gamma} \left|\hat f\big( u(s)\big)\right|^2.$$

\parindent=0mm \vspace{.1in}
which implies,

$$(1-\varepsilon) \big\|f\big\|_{2}^{2} \le \sum_{\lambda =0}^{(qN)^j -1}\left|\left\langle f, \varphi_{j,\lambda}^{\text{per}}\right\rangle\right|^2\le (1+\varepsilon) \big\|f\big\|_{2}^{2}.$$

\parindent=0mm \vspace{.1in}
This completes the proof of Lemma 3.1. ~\fbox

\parindent=0mm \vspace{.2in}

{\bf{Lemma 3.2.}} Let $\varphi$ defined by (2.11) be the nonuniform refinable function with $m_{0}(\xi)$ as its refinement mask and let $m_{\ell}(\xi),  1 \le\ell \le qN-1$ be the wavelet masks. The  nonuniform wavelet system ${\cal W}({\Psi},\lambda)$ given by (2.7) form a normalized tight frame for $L^2(\mathbb K)$. Then, for any function $f \in L^2(\mathbb K)$, we have
$$\sum_{\lambda \in \Lambda}\left|\big\langle f, \varphi_{j+1,\lambda}\big\rangle\right|^2=\sum_{\lambda \in \Lambda}\left|\big\langle f, \varphi_{j,\lambda}\big\rangle\right|^2+\sum_{\ell=1}^{qN-1}\sum_{\lambda \in \Lambda}\left|\big\langle f, \psi_{\ell,j,\lambda}\big\rangle\right|^2.\eqno(3.5)$$

\parindent=0mm \vspace{.1in}
{\bf{Proof.}} For any $f \in L^2(\mathbb K)$ and $j\in \mathbb N_{0}$, define the linear operators ${\cal P}_j$ and ${\cal Q}_j$  as:
$${\cal P}_{j}f(x)=\sum_{\lambda \in \Lambda}\big\langle f, \varphi_{j,\lambda}\big\rangle \varphi_{j,\lambda}(x),\quad {\cal }{\cal Q}_{j}f(x)=\sum_{\ell=1}^{qN-1}\sum_{\lambda \in \Lambda}\big\langle f, \psi_{\ell,j,\lambda}\big\rangle \psi_{\ell,j,\lambda}(x).\eqno(3.6)$$

\parindent=0mm \vspace{.1in}
Since $\Omega(K)$ is dense in $L^2(\mathbb K)$ and closed under Fourier transform, therefore, it is sufficient to prove that
$$\big\langle {\cal P}_jf, f\big\rangle+\big\langle {\cal Q}_jf, f\big\rangle=\big\langle {\cal P}_{j+1}f, f\big\rangle, \eqno(3.7)$$

\parindent=0mm \vspace{.1in}
holds for all the functions $f$ in $\Omega(K)$. Therefore, for all $f\in \Omega(K)$ and $j\in\mathbb Z, k\in\mathbb N_{0}$, we obtain the following equality by using the Parseval's identity
\begin{align*}
\big\langle {\cal P}_jf, f\big\rangle&=(qN)^j\int_{\mathfrak D}\left|\sum_{r \in \mathbb N_0}\hat f\Big(({\mathfrak p}^{-1}N)^{-j}\big(\xi+u(r)\big)\Big)\overline{\hat \varphi\big(\xi+u(r)\big)}\right|^2d \xi\\\\
&=\int_{({\mathfrak p}^{-1}N)^{-j}\mathfrak D}\left|\sum_{r \in \mathbb N_0}\hat f\Big(\xi+({\mathfrak p}^{-1}N)^{-j}u(r)\Big)\overline{\hat \varphi\big(({\mathfrak p}^{-1}N)^{j}\xi+u(r)\big)}\right|^2d \xi.\tag{3.8}
\end{align*}

\parindent=0mm \vspace{.1in}
Since  $m_0(\xi)$ is an integral-periodic function,  equation (3.8) yields
$$\begin{array}{rcl}
\big\langle {\cal P}_jf, f\big\rangle &=&\displaystyle\int_{({\mathfrak p}^{-1}N)^{-j}\mathfrak D}\left|\sum_{r \in \mathbb N_0}\hat f\big(\xi+({\mathfrak p}^{-1}N)^{-j}u(r)\big)\overline{\hat \varphi\big(({\mathfrak p}^{-1}N)^{j+1}\xi+({\mathfrak p}^{-1}N)u(r)\big)}\right.\\\\
&&\qquad\qquad\qquad\qquad\qquad\qquad\qquad \left. \times~ \overline{m_0\big(({\mathfrak p}^{-1}N)^{j+1}\xi+({\mathfrak p}^{-1}N)u(r)\big)}\right|^2d \xi\\\\
&=&\displaystyle\int_{({\mathfrak p}^{-1}N)^{-j}\mathfrak D}\left|\sum_{r \in \mathbb N_0}\sum_{s=0}^{qN-1}\hat f\Big(\xi+({\mathfrak p}^{-1}N)^{-j}\big(({\mathfrak p}^{-1}N)u(r)+u(s)\big)\Big)\right.\\\\
&&\qquad\qquad \left. \times~\overline{\hat \varphi\Big(({\mathfrak p}^{-1}N)^{j+1}\xi+({\mathfrak p}^{-1}N)\big(({\mathfrak p}^{-1}N)u(r)+u(s)\big)\Big)}\right.\\\\
&&\qquad\qquad\qquad\quad \left.\times~\overline{m_0\Big(({\mathfrak p}^{-1}N)^{j+1}\xi+({\mathfrak p}^{-1}N)\big(({\mathfrak p}^{-1}N)u(r)+u(s)\big)\Big)}\right|^2d \xi\\\\
&=&\displaystyle\int_{({\mathfrak p}^{-1}N)^{-j}\mathfrak D}\left|\sum_{r \in \mathbb N_0}\sum_{s=0}^{qN-1}\hat f\Big(\xi+({\mathfrak p}^{-1}N)^{-j}\big(({\mathfrak p}^{-1}N)u(r)+u(s)\big)\Big)\right.\\\\
&&\quad\qquad\qquad\left.\times\overline{\hat \varphi\Big(({\mathfrak p}^{-1}N)^{j+1}\xi+({\mathfrak p}^{-1}N)\big(({\mathfrak p}^{-1}N)u(r)+u(s)\big)\Big)}\right.\\\\
&&\qquad\qquad\qquad\qquad\qquad\qquad\qquad\quad\left.\times~ \overline{m_0\big(({\mathfrak p}^{-1}N)^{j+1}\xi+({\mathfrak p}^{-1}N)u(s)\big)}\right|^2d \xi\\\\
&=&\displaystyle\int_{({\mathfrak p}^{-1}N)^{-j}\mathfrak D}\left|\sum_{s=0}^{qN-1}{\cal R}_{f, \varphi}^j\big(u(s),\xi\big)\overline{m_0\big({\mathfrak p}^{j+1}\xi+{\mathfrak p}u(s)\big)}\right|^2d \xi,
\end{array}$$

\parindent=0mm \vspace{.1in}
where
$$\begin{array}{rcl}
{\cal R}_{f, \varphi}^j\big(u(s),\xi\big)&=&\displaystyle\sum_{r \in \mathbb N_0}\hat f\Big(\xi+({\mathfrak p}^{-1}N)^{-j}\big(({\mathfrak p}^{-1}N)u(r)+u(s)\big)\Big)\\\
&&\qquad\qquad\qquad\times~\overline{\hat \varphi\Big(({\mathfrak p}^{-1}N)^{j+1}\xi+({\mathfrak p}^{-1}N)\big(({\mathfrak p}^{-1}N)u(r)+u(s)\big)\Big)}.
\end{array}$$

\parindent=0mm \vspace{.1in}
Proceeding in the  similar manner as above, we can obtain
$$\big\langle {\cal Q}_{j}f, f\big\rangle=\sum_{\ell=1}^{qN-1}\int_{({\mathfrak p}^{-1}N)^{-j}\mathfrak D}\left|\sum_{s=0}^{qN-1}{\cal R}_{f, \varphi}^j\big(u(s),\xi\big)\overline{m_{\ell}\big(({\mathfrak p}^{-1}N)^{j+1}\xi+({\mathfrak p}^{-1}N)u(s)\big)}\right|^2d \xi.$$

\parindent=0mm \vspace{.1in}
Therefore, we have

\parindent=0mm \vspace{.2in}
$\big\langle {\cal P}_jf, f\big\rangle+\big\langle {\cal Q}_jf, f\big\rangle$
$$\begin{array}{rcl}
&=&\displaystyle\int_{({\mathfrak p}^{-1}N)^{-j}\mathfrak D}\left\{\sum_{s=0}^{qN-1}{\cal R}_{f, \varphi}^j\big(u(s),\xi\big)\overline{m_0\big(({\mathfrak p}^{-1}N)^{j+1}\xi+({\mathfrak p}^{-1}N)u(s)\big)}\right\}\\\\
&& \qquad\qquad\qquad\qquad\times~\overline{\left\{\displaystyle\sum_{s^\prime =0}^{qN-1}{\cal R}_{f, \varphi}^j\big(u(s^\prime),\xi\big)\overline{m_0\big(({\mathfrak p}^{-1}N)^{j+1}\xi+({\mathfrak p}^{-1}N)u(s^\prime)\big)}\right\}}d\xi\\\\
&&\quad+\displaystyle\sum_{\ell=1}^{qN-1}\int_{({\mathfrak p}^{-1}N)^{-j}\mathfrak D}\left\{\sum_{s=0}^{qN-1}{\cal R}_{f, \varphi}^j\big(u(s),\xi\big)\overline{m_{\ell}\big(({\mathfrak p}^{-1}N)^{j+1}\xi+({\mathfrak p}^{-1}N) u(s)\big)}\right\}\\\\
&& \qquad\qquad\qquad\qquad\times~\overline{\left\{\displaystyle\sum_{s^\prime =0}^{qN-1}{\cal R}_{f, \varphi}^j\big(u(s^\prime),\xi\big)\overline{m_{\ell}\big(({\mathfrak p}^{-1}N)^{j+1}\xi+({\mathfrak p}^{-1}N)u(s^\prime)\big)}\right\}}d\xi\\\\
&=&\displaystyle\int_{({\mathfrak p}^{-1}N)^{-j}\mathfrak D}\left\{\sum_{s=0}^{qN-1}\sum_{s^\prime =0}^{qN-1} {\cal R}_{f, \varphi}^j\big(u(s),\xi\big)\overline{{\cal R}_{f, \varphi}^j\big(u(s^\prime),\xi\big)}\right\}\\\\
&& \qquad\times~ \left\{\displaystyle\sum_{\ell=0}^{qN-1} m_{\ell}\big(({\mathfrak p}^{-1}N)^{j+1}\xi+({\mathfrak p}^{-1}N)u(s^\prime)\big)\overline{m_{\ell}\big(({\mathfrak p}^{-1}N)^{j+1}\xi+({\mathfrak p}^{-1}N)u(s)\big)}\right\}d\xi.
\end{array}$$

\parindent=0mm \vspace{.1in}
Since the unitary extension principle condition (2.17) is equivalent to

$$\sum_{\ell=0}^{qN-1} m_{\ell}\left(({\mathfrak p}^{-1}N)^{j+1}\xi+({\mathfrak p}^{-1}N)u(s^\prime)\right)\overline{m_{\ell}\big(({\mathfrak p}^{-1}N)^{j+1}\xi+({\mathfrak p}^{-1}N)u(s)\big)}=\delta_{s,s^\prime}.$$

\parindent=0mm \vspace{.1in}
Therefore, we have

\parindent=0mm \vspace{.2in}
$\big\langle {\cal P}_jf, f\big\rangle+\big\langle {\cal Q}_jf, f\big\rangle$
$$\begin{array}{rcl}
&=&\displaystyle\int_{({\mathfrak p}^{-1})^{-j}\mathfrak D}\sum_{s=0}^{qN-1}\Big|{\cal R}_{f, \varphi}^j\big(u(s),\xi\big)\Big|^2d \xi\\\\
&=&\displaystyle\int_{({\mathfrak p}^{-1})^{-j}\mathfrak D}
\sum_{s=0}^{qN-1}\left|\sum_{r \in \mathbb N_0}\hat f\Big(\xi+({\mathfrak p}^{-1}N)^{-j}\big(({\mathfrak p}^{-1}N)u(r)+u(s)\big)\Big)\right.\\\

&&\qquad\qquad\qquad\qquad\left.\times~\overline{\hat \varphi\Big(({\mathfrak p}^{-1}N)^{j+1}\xi+({\mathfrak p}^{-1}N)\big(({\mathfrak p}^{-1}N)u(r)+u(s)\big)\Big)}\right|^2d \xi\\\\
&=&\displaystyle\sum_{s=0}^{qN-1}\int_{({\mathfrak p}^{-1}N)^{-j}{\mathfrak D}+({\mathfrak p}^{-1}N)^{-j}u(s)}\left|\displaystyle\sum_{r \in \mathbb N_0}\hat f\big(\xi+({\mathfrak p}^{-1}N)^{-j-1}u(r)\big)\overline{\hat \varphi\big(({\mathfrak p}^{-1}N)^{j+1}\xi+u(r)\big)}\right|^2d \xi\\\\
&=&\displaystyle\int_{({\mathfrak p}^{-1}N)^{-j-1}{\mathfrak D}}\left|\displaystyle\sum_{r \in \mathbb N_0}\hat f\big(\xi+({\mathfrak p}^{-1}N)^{-j-1}u(r)\big)\overline{\hat \varphi\big(({\mathfrak p}^{-1}N)^{j+1}\xi+u(r)\big)}\right|^2d \xi\\\\
&=&\big\langle {\cal P}_{j+1}f, f\big\rangle,
\end{array}$$

\parindent=0mm \vspace{.1in}
and hence, we get the desired result (3.5). \qquad \fbox

\parindent=0mm \vspace{.2in}

{\bf{Lemma 3.3.}}  Let $\varphi\in L^2(\mathbb K)$ be a nonuniform refinable function with refinement mask $m_0(\xi)$, and let the wavelet system ${\cal W}({\Psi},\lambda)$ given by (2.7) constitutes a normalized tight frame for $L^2(\mathbb K)$. If  $ \left\{\varphi, \psi_1, \dots,\psi_L\right\}\subset L^1(\mathbb K)\cap L^2(\mathbb K)$ and $\varphi, \psi_1, \dots,\psi_{qN-1}$ have a common radial decreasing $L^{1}$-majorant, then we have
$$\sum_{\lambda=0}^{(qN)^{j+1}-1}\left|\left\langle f, \varphi_{j+1,\lambda}^{\text{per}}\right\rangle\right|^2=\sum_{\lambda=0}^{(qN)^j -1}\left|\left\langle f, \varphi_{j,\lambda}^{\text{per}}\right\rangle\right|^2+\sum_{\ell=1}^{qN-1}\sum_{\lambda =0}^{(qN)^{j} -1}\left|\left\langle f, \psi_{\ell,j,\lambda}^{\text{per}}\right\rangle\right|^2.\eqno(3.9)$$  

\parindent=0mm \vspace{.0in}
{\bf{Proof.}} For any $f\in \Omega(K)$ and $j \in \mathbb N_{0}$, we have
$$\sum_{\lambda =0}^{(qN)^j-1}\left|\left\langle f, \varphi_{j,\lambda}^{\text{per}}\right\rangle\right|^2=\sum_{\lambda  =0}^{(qN)^j-1}\left|\Big\langle f, \sum_{r \in \mathbb N_0}\varphi_{j,\lambda}\big(x+u(r)\big)\Big\rangle\right|^2=\sum_{\lambda =0}^{(qN)j-1}\left|\sum_{r\in \mathbb N_0}\Big\langle f, \varphi_{j,\lambda}\big(x+u(r)\big)\Big\rangle\right|^2.$$

\parindent=0mm \vspace{.1in}
Furthermore, we have
$$\begin{array}{rcl}
\displaystyle\sum_{r\in \mathbb N_0}\int_{\mathfrak D}\left|f(x) \overline{\varphi_{j,\lambda}\big(x+u(r)\big)}\right|dx &\le& \displaystyle \left\|f\right\|_{L^{\infty}(\mathfrak D)}\int_\mathbb K\big|\varphi_{j,\lambda}(x)\big|dx\\\
&=&\displaystyle \left|\big|f\big|\right|_{L^{\infty}(\mathfrak D)}(qN)^{j/2}\int_\mathbb K\big|\varphi(x)\big|dx<\infty.
\end{array}$$

\parindent=0mm \vspace{.0in}
Which implies that  the series

$$\sum_{\lambda =0}^{(qN^j -1 }\sum_{r\in \mathbb N_0}\sum_{s\in \mathbb N_0}\Big\langle f, \varphi_{j,\lambda}\big(x+u(r)\big)\Big\rangle\overline{\Big\langle f, \varphi_{j,\lambda}\big(x+u(s)\big)\Big\rangle}$$

\parindent=0mm \vspace{.1in}
is absolutely convergent. Therefore, the series can be rearranged as follows:
$$\begin{array}{rcl}
\displaystyle\sum_{\lambda =0}^{(qN)^j -1}\left|\left\langle f, \varphi_{j,\lambda}^{\text{per}}\right\rangle\right|^2&=&\displaystyle\sum_{\lambda =0}^{(qN)^j -1}\sum_{r\in \mathbb N_0}\sum_{s\in \mathbb N_0}\Big\langle f, \varphi_{j,\lambda}\big(x+u(r)\big)\Big\rangle\overline{\Big\langle f, \varphi_{j,\lambda}\big(x+u(s)\big)\Big\rangle}\\\\
&=&\displaystyle\sum_{\lambda=0}^{(qN)^j -1}\sum_{r\in \mathbb N_0}\sum_{s\in \mathbb N_0}\Big\langle f, \varphi_{j,\lambda}\big(x+u(r)\big)\Big\rangle\overline{\Big\langle f, \varphi_{j,\lambda}\big(x+u(r)+u(s)\big)\Big\rangle}.
\end{array}$$

\parindent=0mm \vspace{.1in}
For $s\in \mathbb N_0$, we define
$$F_s(x)=f(x)\,{\bf 1}_{{\mathfrak D}+u(s)}(x),$$

\parindent=0mm \vspace{.1in}
where ${\bf 1}_{{\mathfrak D}+u(s)}(x)$ is the characteristic function of the set ${\mathfrak D}+u(s)$. Using the fact that $\varphi_{j,\lambda}\big(x+u(s)\big)=\varphi_{j,\lambda-({\mathfrak p}^{-1}N)^{-j}u(s)}(x),$ we have
$$\begin{array}{rcl}
\displaystyle\sum_{\lambda =0}^{(qN)^j -1}\left|\left\langle f, \varphi_{j,\lambda}^{\text{per}}\right\rangle\right|^2 &=&\displaystyle\sum_{\lambda =0}^{(qN)^j -1}\sum_{r\in \mathbb N_0}\sum_{s\in \mathbb N_0}\left\{\int_{\mathfrak D} f(x)\overline{ \varphi_{j,\lambda}\big(x+u(r)\big)}dx\right\}\\\
&&\qquad\qquad\qquad\qquad\qquad \times~\overline{\left\{\displaystyle\int_{\mathfrak D} f(x)\overline{ \varphi_{j,\lambda}\big(x+u(r)+u(s)\big)}dx\right\}}\\\\
&=&\displaystyle\sum_{\lambda =0}^{(qN)^j -1}\sum_{r\in \mathbb N_0}\sum_{s\in \mathbb N_0}\left\{\int_{\mathfrak D} f(x)\overline{ \varphi_{j,\lambda}\big(x+u(r)\big)}dx\right\}\\\
&&\qquad\qquad\qquad\qquad\qquad\qquad \times~\overline{\left\{\displaystyle\int_{\mathfrak D+u(s)} f(x)\overline{ \varphi_{j,\lambda}\big(x+u(r)\big)}dx\right\}}\\\\
&=&\displaystyle\sum_{\lambda  =0}^{(qN)^j -1}\sum_{r\in \mathbb N_0}\sum_{s\in \mathbb N_0}\left\{\int_{\mathbb K} F_0(x)\overline{ \varphi_{j,\lambda}\big(x+u(r)\big)}dx\right\}\\\
&&\qquad\qquad\qquad\qquad\qquad\qquad\qquad \times~\overline{\left\{\displaystyle\int_{\mathbb K} F_s(x)\overline{ \varphi_{j,\lambda}\big(x+u(r)\big)}dx\right\}}\\\\
&=&\displaystyle\sum_{\lambda =0}^{(qN)^j-1}\sum_{r\in \mathbb N_0}\sum_{s\in \mathbb N_0}\left\langle F_0, \varphi_{j,\lambda-({\mathfrak p}^{-1}N)^{-j}u(r)}\right\rangle\overline{\left\langle F_s, \varphi_{j,\lambda-({\mathfrak p}^{-1}N)^{-j}u(r)}\right\rangle}\\\\
&=&\displaystyle\sum_{\lambda\in \Lambda}\sum_{s\in \mathbb N_0}\left\langle F_0, \varphi_{j,\lambda}\right\rangle\overline{\left\langle F_s, \varphi_{j,\lambda}\right\rangle}.
\end{array}$$

\parindent=0mm \vspace{.1in}
Similarly, for each $ 1 \le \ell \le qN-1$, we have

$$\sum_{\lambda =0}^{(qN)^j-1}\left|\left\langle f, \psi_{\ell,j,\lambda}^{\text{per}}\right\rangle\right|^2=\sum_{\lambda\in \Lambda}\sum_{s\in \mathbb N_0}\left\langle F_0, \psi_{\ell,j,\lambda}\right\rangle\overline{\left\langle F_s, \psi_{\ell,j,\lambda}\right\rangle}.$$

\parindent=0mm \vspace{.1in}
By Lemma 3.2, we have
$$\begin{array}{lcr}
\displaystyle\sum_{\lambda=0}^{(qN)^j-1}\left|\left\langle f, \varphi_{j,\lambda}^{\text{per}}\right\rangle\right|^2+\sum_{\ell=1}^{qN-1}\sum_{\lambda =0}^{(qN)^j-1}\left|\left\langle f, \psi_{\ell,j,\lambda}^{\text{per}}\right\rangle\right|^2&&\\\
\qquad\qquad\qquad=\displaystyle\sum_{\lambda\in \Lambda}\sum_{s\in \mathbb N_0}\left\langle F_0, \varphi_{j,\lambda}\right\rangle\overline{\left\langle F_s, \varphi_{j,\lambda}\right\rangle}+\sum_{\ell=1}^{qN-1}\sum_{\lambda \in \Lambda}\sum_{s\in \mathbb N_0}\left\langle F_0, \psi_{\ell,j,\lambda}\right\rangle\overline{\left\langle F_s, \psi_{\ell,j,\lambda}\right\rangle}&&\\\\
\qquad\qquad\qquad=\displaystyle\sum_{\lambda \in \Lambda}\sum_{s\in \mathbb N_0}\left\langle F_0, \varphi_{j+1,\lambda}\right\rangle\overline{\left\langle F_s, \varphi_{j+1,\lambda}\right\rangle}&&\\\\
\qquad\qquad\qquad=\displaystyle\sum_{\lambda=0}^{(qN)^j -1}\left|\left\langle f, \varphi_{j+1,\lambda}^{\text{per}}\right\rangle\right|^2.
\end{array}$$

\parindent=0mm \vspace{.1in}
This completes the proof of Lemma 3.2. \qquad \fbox

\parindent=0mm \vspace{.1in}

Now we state and prove the  main result of this section.

\parindent=0mm \vspace{.2in}
{\bf{Theorem 3.1.}}  Let $m_{0}(\xi)$ be the refinement mask of a refinable function $\varphi(x)$  and let $m_{\ell}(\xi),  1 \le \ell \le qN-1$ be the wavelet masks associated with the basic wavelets given by (2.13). Furthermore, let the wavelet system ${\cal W}({\Psi},\lambda)$ given by (2.7) form a normalized tight
frame generated by the refinable function $\phi$. If  $ \left\{\varphi, \psi_1, \psi_2,\dots,\psi_{qN-1}\right\}\subset L^1(K)\cap L^2(K)$ and $\varphi, \psi_1, \psi_2,\dots,\psi_{qN-1}$ have a common radial decreasing $L^{1}$-majorant, then the periodic wavelet system ${\cal W}^{\text{per}}({\Psi},\lambda)$ given by (3.1) generates a normalized tight frame for $L^2(\mathfrak D)$. 

\parindent=0mm \vspace{.1in}
{\bf{Proof.}} For any periodic function $f \in \Omega(\mathfrak D)$ and $\varepsilon>0$, we can choose $J>0$ by Lemma 3.1 such that for all $j>J$, we have
$$(1-\varepsilon)\big\|f\big\|_2^2 \le \sum_{\lambda=0}^{(qN)^j-1}\left|\left\langle f, \varphi_{j,k}^{\text{per}}\right\rangle\right|^2 \le  (1+\varepsilon)\big\|f\big\|_2^2.$$

\parindent=0mm \vspace{.1in}
Also for any $j \in \mathbb Z$, Lemma 3.3 implies that

$$\sum_{\lambda =0}^{(qN)^j-1}\left|\left\langle f, \varphi_{j,\lambda}^{\text{per}}\right\rangle\right|^2=\sum_{\lambda=0}^{(qN){j-1}-1}\left|\left\langle f, \varphi_{j-1,\lambda}^{\text{per}}\right\rangle\right|^2+\sum_{\ell=1}^{qN-1}\sum_{\lambda=0}^{(qN)^{j-1}-1}\left|\left\langle f, \psi_{\ell,j-1,\lambda}^{\text{per}}\right\rangle\right|^2.$$

\parindent=0mm \vspace{.1in}
Allowing the above argument to repeat on $\displaystyle\sum_{\lambda=0}^{(qN)^{j-1}-1}\left|\left\langle f, \varphi_{j-1,\lambda}^{\text{per}}\right\rangle\right|^2$, we obtain

$$\sum_{\lambda =0}^{(qN)^j -1}\left|\left\langle f, \varphi_{j,\lambda}^{\text{per}}\right\rangle\right|^2=\big|\big\langle f, \varphi^{\text{per}}\big\rangle\big|^2+\sum_{\ell=1}^{qN-1}\sum_{r=0}^{j-1}\sum_{\lambda=0}^{(qN)^r -1}\left|\left\langle f, \psi_{\ell,r,\lambda}^{\text{per}}\right\rangle\right|^2.$$

\parindent=0mm \vspace{.1in}
Therefore, we have

$$(1-\varepsilon)\big\|f\big\|_2^2 \le \big|\big\langle f, \varphi^{\text{per}}\big\rangle\big|^2+\sum_{\ell=1}^{qN-1}\sum_{r=0}^{j-1}\sum_{\lambda=0}^{(qN)^r -1}\left|\left\langle f, \psi_{\ell,r,\lambda}^{\text{per}}\right\rangle\right|^2 \le  (1+\varepsilon)\big\|f\big\|_2^2.$$

\parindent=0mm \vspace{.1in}
Letting $j \to \infty$, we obtain

$$(1-\varepsilon)\big\|f\big\|_2^2 \le \big|\big\langle f, \varphi^{\text{per}}\big\rangle\big|^2+\sum_{\ell=1}^{qN-1}\sum_{r\in \mathbb N_{0}}\sum_{\lambda =0}^{(qN)^r-1}\left|\left\langle f, \psi_{\ell,r,\lambda}^{\text{per}}\right\rangle\right|^2 \le  (1+\varepsilon)\big\|f\big\|_2^2.$$

\parindent=0mm \vspace{.1in}
Since $\epsilon>0$ was arbitrary. Hence, it follows that

$$\big|\big\langle f, \varphi^{\text{per}}\big\rangle\big|^2+\sum_{\ell=1}^{qN-1}\sum_{r\in \mathbb N_{0}}\sum_{\lambda=0}^{(qN)^r-1}\left|\left\langle f, \psi_{\ell,r,\lambda}^{\text{per}}\right\rangle\right|^2=\|f\big\|_2^2.$$

\parindent=0mm \vspace{.1in}
This completes the proof of  Theorem 3.4. \qquad \fbox

\parindent=0mm \vspace{.2in}

{\bf{References}}


\begin{enumerate}


{\small {

\bibitem{onuwf} O. Ahmad and N. A. Sheikh,  On Characterization of nonuniform tight  wavelet frames on local fields,{\it  Anal. Theory Appl.,}  {\bf 34} (2018) 135-146.

\bibitem{BB} J.J. Benedetto and R.L. Benedetto, A wavelet theory for local fields and related groups. {\it J. Geom. Anal.} {\bf 14} (2004) 423-456.

\bibitem{chrisBOOK} O. Christensen,  {\it An Introduction to Frames and Riesz Bases}, Second Edition, Birkh\"{a}user, Boston, 2016.

\bibitem{DHRS} I. Daubechies, B. Han, A. Ron and Z. Shen,  Framelets: MRA-based constructions of wavelet frames, {\it Appl. Comput. Harmon. Anal.} {\bf 14} (2003) 1-46.

\bibitem{DGM} I. Daubechies, A. Grossmann, Y. Meyer,  Painless non-orthogonal expansions,{\it J. Math. Phys.} {\bf 27} (5) (1986) 1271-1283.

\bibitem{DS} R.J. Duffin and A.C. Shaeffer, A class of nonharmonic Fourier series. {\it Trans. Amer. Math. Soc.} {\bf 72} (1952) 341-366.

\bibitem{GN1} J. P. Gabardo and M. Nashed,  Nonuniform multiresolution analyses and spectral pairs, {\it J. Funct. Anal.} {\bf 158} (1998)  209-241.

\bibitem{GN2} J.P. Gabardo and X. Yu, Wavelets associated with nonuniform multiresolution analysis and one-dimensional spectral pairs, {\it J. Math. Anal. Appl.}  323 (2006), 798-817.

\bibitem{jiang} H.K. Jiang, D.F. Li and N. Jin,  Multiresolution analysis on local fields. {\it J. Math. Anal. Appl.} {\bf 294} (2004) 523-532.

\bibitem{LP} E. A. Lebedeva and J. Prestin,  Periodic wavelet frames and time-frequency localization. {\it Appl. Comput. Harmon. Anal.} {\bf 37} (2014) 347-359.

\bibitem{LL} D. Lu and D. Li,  Construction of periodic wavelet frames with dilation matrix. {\it Front. Math. China.} {\bf 9} (2014) 111-134

\bibitem{RV} D. Ramakrishnan and R. J. Valenza,  Fourier Analysis on Number Fields, Graduate Texts in Mathematics 186, Springer-Verlag, New York (1999).

\bibitem{RS} A. Ron and Z. Shen,  Affine systems in $L^2(\mathbb R^d)$: the analysis of the analysis operator. {\it J.  Funct. Anal.} {\bf 148} (1997) 408-447.

\bibitem{oJMP} F. A. Shah, O. Ahmad and P.E. Jorgenson, Fractional Wave Packet Frames in  $L^2(\mathbb R)$,{\it Journal of Mathematical Physics,} 59, 073509 (2018); doi: 10.1063/1.5047649.

\bibitem{oJGP} F.A. Shah  and O.  Ahmad,  Wave packet systems on local fields,{\it Journal of Geometry and Physics},  {\bf 120}  (2017) 5-18.

\bibitem{SDtight} F.A. Shah and L. Debnath. Tight wavelet frames on local fields. {\it Analysis}.  {\bf 33} (2013) 293-307.


\bibitem{table} M.H. Taibleson,  {\it Fourier Analysis on Local Fields}, Princeton University Press, Princeton, NJ,(1975).

\bibitem{zhang} Z. Zhang, Periodic wavelet frames. {\it Adv. Comput. Math.} {\bf 22} (2005) 165-180.

\bibitem{ZS} Z. Zhang and N. Saito,  Constructions of periodic wavelet frames using extension principles. {\it Appl. Comput. Harmon. Anal.} {\bf 27} (2009) 12-23.


}}
\end{enumerate}
\end{document}